\documentclass[11pt]{amsart}
\usepackage{amsmath,amsthm,epsfig,amssymb}
\usepackage{xypic}
\input xy
\title{Counting all regular tetrahedra in $\{0,1,...,n\}^3$}
\author{Eugen J. Ionascu}
\curraddr{Department of Mathematics\\ Columbus State University\\4225 University Avenue\\
Columbus, GA 31907\\
Honorific Member of the Romanian Institute of Mathematics ``Simion
Stoilow" } \email{ionascu\_eugen@colstate.edu;} \subjclass{}
\date{December $5^{th}$, 2009}
\begin{document}
\def\sms{\small\scshape}
\baselineskip18pt
\newtheorem{theorem}{\hspace{\parindent}
T{\scriptsize HEOREM}}[section]
\newtheorem{proposition}[theorem]
{\hspace{\parindent }P{\scriptsize ROPOSITION}}
\newtheorem{corollary}[theorem]
{\hspace{\parindent }C{\scriptsize OROLLARY}}
\newtheorem{lemma}[theorem]
{\hspace{\parindent }L{\scriptsize EMMA}}
\newtheorem{definition}[theorem]
{\hspace{\parindent }D{\scriptsize EFINITION}}
\newtheorem{problem}[theorem]
{\hspace{\parindent }P{\scriptsize ROBLEM}}
\newtheorem{conjecture}[theorem]
{\hspace{\parindent }C{\scriptsize ONJECTURE}}
\newtheorem{example}[theorem]
{\hspace{\parindent }E{\scriptsize XAMPLE}}
\newtheorem{remark}[theorem]
{\hspace{\parindent }R{\scriptsize EMARK}}
\renewcommand{\thetheorem}{\arabic{section}.\arabic{theorem}}
\renewcommand{\theenumi}{(\roman{enumi})}
\renewcommand{\labelenumi}{\theenumi}
\newcommand{\Q}{{\mathbb Q}}
\newcommand{\Z}{{\mathbb Z}}
\newcommand{\N}{{\mathbb N}}
\newcommand{\C}{{\mathbb C}}
\newcommand{\R}{{\mathbb R}}
\newcommand{\F}{{\mathbb F}}
\newcommand{\K}{{\mathbb K}}
\newcommand{\D}{{\mathbb D}}
\def\phi{\varphi}
\def\ra{\rightarrow}
\def\sd{\bigtriangledown}
\def\ac{\mathaccent94}
\def\wi{\sim}
\def\wt{\widetilde}
\def\bb#1{{\Bbb#1}}
\def\bs{\backslash}
\def\cal{\mathcal}
\def\ca#1{{\cal#1}}
\def\Bbb#1{\bf#1}
\def\blacksquare{{\ \vrule height7pt width7pt depth0pt}}
\def\bsq{\blacksquare}
\def\proof{\hspace{\parindent}{P{\scriptsize ROOF}}}
\def\pofthe{P{\scriptsize ROOF OF}
T{\scriptsize HEOREM}\  }
\def\pofle{\hspace{\parindent}P{\scriptsize ROOF OF}
L{\scriptsize EMMA}\  }
\def\pofcor{\hspace{\parindent}P{\scriptsize ROOF OF}
C{\scriptsize ROLLARY}\  }
\def\pofpro{\hspace{\parindent}P{\scriptsize ROOF OF}
P{\scriptsize ROPOSITION}\  }
\def\n{\noindent}
\def\wh{\widehat}
\def\eproof{$\hfill\bsq$\par}
\def\ds{\displaystyle}
\def\du{\overset{\text {\bf .}}{\cup}}
\def\Du{\overset{\text {\bf .}}{\bigcup}}
\def\b{$\blacklozenge$}

\def\eqtr{{\cal E}{\cal T}(\Z) }
\def\eproofi{\bsq}

\begin{abstract} In this note we describe a procedure of calculating the number all regular
tetrahedra that have coordinates in the set $\{0,1,...,n\}$. We
develop a few results that may help in finding good estimates for
this sequence which is twice A103158 in the Online Encyclopedia of
Integer Sequences \cite{OL}.
\end{abstract} \maketitle
\section{INTRODUCTION}
The story of regular tetrahedra having vertices of integer
coordinates starts with the parametrization of some equilateral
triangles in $\mathbb Z^3$ that begun in \cite{eji}. There was an
additional hypothesis that did not cover all the generality in the
result obtained in \cite{eji} but it was removed successfully in
\cite{rceji}. A few other related results appeared in \cite{ejic}
and \cite{ejirt}. In this note we are interested in the following
problem.

{\it How many regular tetrahedra, $T(n)$, can be found if the
coordinates of its vertices must be in the set $\{0,1,...,n\}$? We
observe that $A103158=\frac{1}{2}T(n)$, see \cite{OL}.} This
sequence starts as in the next tables:

\vspace{0.5in}

\centerline{ \vspace{0.2in}
\begin{tabular}{|c||c|c|c|c|c|c|c|c|c|c|c|}
  \hline
  n& 1 & 2 & 3 & 4 & 5 & 6 & 7 & 8 & 9 &10&11\\ \hline
 A103158& 1 & 9 & 36 & 104 & 257 & 549 & 1058 & 1896 & 3199&5154 &7926  \\
  \hline
\end{tabular}}
\vspace{0.1in}

\centerline{
\begin{tabular}{|c||c|c|c|c|c|c|c|}
  \hline
  n &12&13& 14& 15 & 16 & 17 & 18 \\ \hline
 A103158&11768 & 16967& 23859& 32846 & 44378 & 58977 & 77215 \\
  \hline
\end{tabular}.}
\vspace{0.1in}

Using our method which is going to be described later we extended
this sequence for all $n\le 100$, and one can go far enough with
this if time allows and powerful computer is used.

\vspace{0.1in}

The rest of the terms are included at the end of the paper. Our
approach begins with looking first at the faces of a regular
tetrahedron, which must be equilateral triangles. It turns out that
every equilateral triangle in $\mathbb Z^3$ after a translation by a
vector with integer coordinates can be assumed to have the origin as
one of its vertices. Then one can show that the other triangle's
vertices are contained in a lattice of points of the form

\begin{equation}\label{planelattice}
{\cal P}_{a,b,c}:=\{(\alpha,\beta,\gamma)\in \mathbb {Z}^3|\ \
a\alpha+b\beta+c\gamma=0,\ \ a^2+b^2+c^2=3d^2,\ \ a,b,c,d\in \mathbb
Z\}.
\end{equation}

\begin{center}\label{figure0}
\epsfig{file=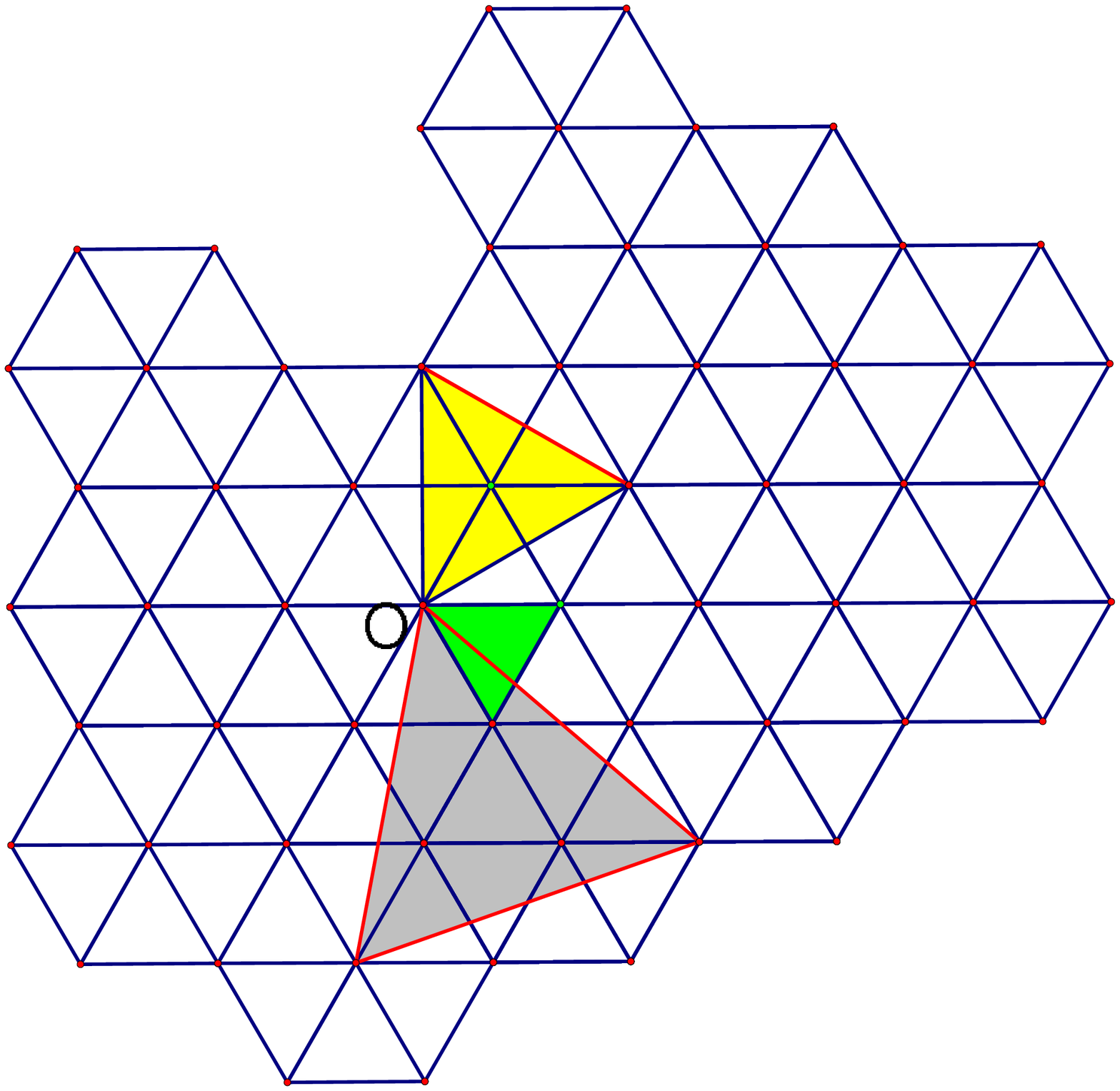,height=3in,width=3in}
\vspace{.1in}\\
{\small Figure 1: The lattice ${\cal P}_{a,b,c}$ } \vspace{.1in}
\end{center}

\n In general, the vertices of the equilateral triangles that dwell
in ${\cal P}_{a,b,c}$ form a strict sub-lattice of ${\cal
P}_{a,b,c}$  which is generated by only two vectors,
$\overrightarrow{\zeta}$ and $\overrightarrow{\eta}$ (see
Figure~\ref{figure0}). These two vectors are described by the
Theorem~\ref{generalpar} proved in \cite{rceji}.

\begin{theorem}\label{generalpar} Let $a$, $b$, $c$, $d$ be odd
integers such that $a^2+b^2+c^2=3d^2$ and ${\rm gcd}(a,b,c)=1$. Then
for every $m,n \in\mathbb{Z}$ (not both zero) the triangle $OPQ$,
determined by

\begin{equation}\label{vectorid}
\overrightarrow{OP}=m\overrightarrow{\zeta}-n\overrightarrow{\eta},\
\
\overrightarrow{OQ}=n\overrightarrow{\zeta}-(n-m)\overrightarrow{\eta},
\ { \rm with} \ \overrightarrow{\zeta}=(\zeta_1,\zeta_1,\zeta_2),
\overrightarrow{\eta}=(\eta_1,\eta_2,\eta_3),
\end{equation}

\begin{equation}\label{paramtwo}
\begin{array}{l}
\begin{cases}
\ds \zeta_1=-\frac{rac+dbs}{q}, \\ \\
\ds \zeta_2=\frac{das-bcr}{q},\\ \\
\ds \zeta_3=r,
\end{cases}
,\ \
\begin{cases}
\ds \eta_1=-\frac{db(s-3r)+ac(r+s)}{2q},\\ \\
\ds \eta_2=\frac{da(s-3r)-bc(r+s)}{2q},\\ \\
\ds \eta_3=\frac{r+s}{2},
\end{cases}
\end{array}
\end{equation}
where $q=a^2+b^2$ and $(r,s)$ is a suitable solution of
$2q=s^2+3r^2$ that makes all the numbers in {\rm (\ref{paramtwo})}
integers, forms an equilateral triangle in $\mathbb Z^3$ contained
in the lattice {\rm (\ref{planelattice})} and having sides-lengths
equal to $d\sqrt{2(m^2-mn+n^2)}$.

Conversely, there exists a choice of the integers $r$ and $s$ such
that given an arbitrary equilateral triangle in $\mathbb{R}^3$ whose
vertices, one at the origin and the other two in the lattice  {\rm
(\ref{planelattice})}, then there also exist integers $m$ and $n$
such that the two vertices not at the origin are given by {\rm
(\ref{vectorid})} and {\rm (\ref{paramtwo})}.
\end{theorem}

The Diophantine equation

\begin{equation}\label{plane}
a^2+b^2+c^2=3d^2
\end{equation}

\n has non-trivial solutions for every $d$ odd. As a curiosity, for
$d=2009$ one obtains 294 solutions satisfying also $0<a\le b\le c$
and  gcd$(a,b,c)=1$. We will refer to such a solution of
(\ref{plane}) as a {\it positive ordered primitive} solution. For
$d=2008$, all these solutions satisfy even a stronger condition:
$a<b<c$. Determining the exact number of solutions for (\ref{plane})
is certainly important if one wishes to find the number (or just an
estimate) of equilateral triangles or the number of tetrahedra with
vertices in $\{0,1,2,...,n\}^3$ . The number of solutions for
(\ref{plane}), coincidentally, taken into account all permutations
and changes of signs is given in a 1999 paper of Hirschhorn and
Seller \cite{hs}:

\begin{equation}\label{numberofsol}
8\left[\underset{
\begin{array}{c}
  p\equiv 1\ or\  7 (mod\ 12) \\
  p^\beta ||d
\end{array}
} {\prod}p^{\beta} \right]\left[\underset{
\begin{array}{c}
  q\equiv 5\ or\ 11 (mod\ 12)\\
  q^\alpha ||d
\end{array}
} {\prod}\left
(q^{\alpha}+2\frac{q^{\alpha}-1}{q-1}\right)f(d)\right],
\end{equation}

\n where $f(d)=\begin{cases} 1 \ if \ 3|d \\
\frac{3^{\gamma}-1}{2}\ if\ 3^{\gamma}||d.  \end{cases}$ Even more
important for our purpose is the calculation of the number of
primitive representations of $d$ as in (\ref{plane})
(gcd$(a,b,c)=1$) in terms of $d$ which appeared in a very recent
paper of Cooper and Hirschhorn \cite{ch}. One may check easily that
the following is  a corollary of Theorem~2 in \cite{ch}.

\begin{theorem}\label{primitivesolch}
{\bf [Cooper-Hirschhorn]} Given an odd number $d$, the number of
primitive solutions of {\rm (\ref{plane})} taking into account all
changing of sings and permutations,  is equal to

\begin{equation}\label{2007hirschhorn}
\Lambda(d):=8d\underset{p|d, p\
prime}{\prod}\left(1-\frac{(\frac{-3}{p})}{p}\right),
\end{equation}

where $(\frac{-3}{p})$ it the Legendre symbol.
\end{theorem}

We remind the reader that, if $p$ is prime then

\begin{equation}\label{legendresymbol}
(\frac{-3}{p})=\begin{cases}0\ \  {\rm if }\  p=3\\ \\
1\ \  {\rm if }\  p\equiv 1\ {\rm or\ } 7 \ {\rm (mod \ 12)} \\ \\
-1\ \  {\rm if }\  p\equiv 5\ {\rm or\ } 11 \ {\rm (mod \ 12)}
 \ \end{cases}.
\end{equation}

\n We observe that the same type of prime partition are used into
different calculations in both formulae (\ref{numberofsol}) and
(\ref{2007hirschhorn}). We have mentioned that the number of
positive ordered primitive representations for $d=2009$ was 294.
This is exactly the number given by (\ref{2007hirschhorn}) modulo
the number of permutations and changes of signs: indeed,
$2009=(41)(7^2)$, $(\frac{-3}{7})=1$, $(\frac{-3}{41})=-1$ and
$\frac{\Lambda(2009)}{48}=\frac{8(41)(7^2)}{48}(1-\frac{1}{7})(1+\frac{1}{41})=7(42)=294.$
This happens because there are no repeating values for $a$, $b$ and
$c$ in any of the positive ordered primitive solutions of
(\ref{plane}). We will see later how the correct number or positive
ordered primitive representations can be obtained in general.

For $k\in \mathbb N$, we let $\Omega:=\{(m,n)\in \mathbb Z\times
\mathbb Z: m^2-mn+n^2=k^2\}$.  In \cite{eji} we showed that every
regular tetrahedron with integer coordinates must have side lengths
of the form  $\lambda\sqrt{2}$, $\lambda\in \mathbb N$, and in
\cite{ejirt} we have found the following characterization of the
regular tetrahedrons with integer coordinates.

\begin{theorem}\label{main} Every tetrahedron whose side lengths are $\lambda\sqrt{2}$,
$\lambda\in \mathbb N$, which has a vertex at the origin,  can be
obtained by taking as one of its faces an equilateral triangle
having the origin as a vertex and the other two vertices given by
{\rm (\ref{vectorid})} and {\rm (\ref{paramtwo})} with $a$, $b$, $c$
and $d$ odd integers satisfying {\rm (\ref{plane})} with $d$ a
divisor of $\lambda$, and then completing it with the fourth vertex
$R$ with coordinates

\begin{equation}\label{fourthpoint}
\begin{array}{c}
 \displaystyle \left(
\frac{\begin{array}{c}
        (2\zeta_1-\eta_1)m \\
        -(\zeta_1+\eta_1)n \\
       \pm 2ak
      \end{array}
}{3},\frac{\begin{array}{c}
        (2\zeta_2-\eta_2)m \\
        -(\zeta_2+\eta_2)n \\
       \pm 2bk
      \end{array}}{3},
 \displaystyle \frac{\begin{array}{c}
        (2\zeta_3-\eta_3)m \\
        -(\zeta_3+\eta_3)n \\
       \pm 2ck
      \end{array}}{3} \right),  \ for\  some\  (m,n)\in
\Omega(k),\ k:=\frac{\lambda }{d}.
\end{array}
\end{equation}

\n Conversely, if we let $a$, $b$, $c$ and $d$ be a primitive
solution of {\rm (\ref{plane})}, let $k\in \Bbb N$ and $(m,n)\in
\Omega(k)$, then the coordinates of the point $R$ in {\rm
(\ref{fourthpoint})}, which completes the equilateral triangle $OPQ$
given as in {\rm (\ref{vectorid})} and {\rm (\ref{paramtwo})}, are

(a) all integers, if $k\equiv 0$ {\rm (mod  3)} regardless of the
choice of signs or

(b) integers, precisely for only one choice of the signs if $k\not
\equiv 0$ {\rm (mod 3)}.
\end{theorem}

The following graph  (Figure 2) is constructed on the positive
ordered primitive solutions of (\ref{plane}), with edges defined by:
\begin{quote}{\em two vertices, say $[(a_1,b_1,c_1),d_1]$ and
$[(a_2,b_2,c_2),d_2]$, are connected, if and only if

\begin{equation}\label{angleBetweenPlanes}
a_1a'_2\pm a_2b'_2\pm c_1c'_2\pm d_1d_2=0\end{equation}

\n for some choice of the signs and permutation $(a'_2,b'_2,c'_2)$
of $(a_2,b_2,c_2)$.}
\end{quote}

Equation (\ref{angleBetweenPlanes}) insures basically that the
planes $\cal P_{a_1,b_1,c_1}$ and $\cal P_{a_2',b_2',c_2'}$
associated to two faces make a dihedral angle of $\arccos
(1/3)\approx 70.52878^{\circ}$. In fact, this equality characterizes
the existence of a regular tetrahedron having integer coordinates
with one of its faces in the plane $\cal P_{a_1,b_1,c_1}$ and
another contained in the plane $\cal P_{a_2',b_2',c_2'}$.  For
instance, $[(1,1,5),3]$ is connected to $[(1,5,11),7]$ since
$1(11)+(1)5+5(1)-3(7)=0$. An example of a regular tetrahedron which
has a face in $\cal P_{-5,-1,1}$ and one face in $\cal P_{-1,-5,11}$
is given by the vertices: $[19, 23, 0]$, $[0, 12, 20]$, $[27, 0,
17]$, and $[24, 27, 29]$.

\begin{center}\label{figure1}
\epsfig{file=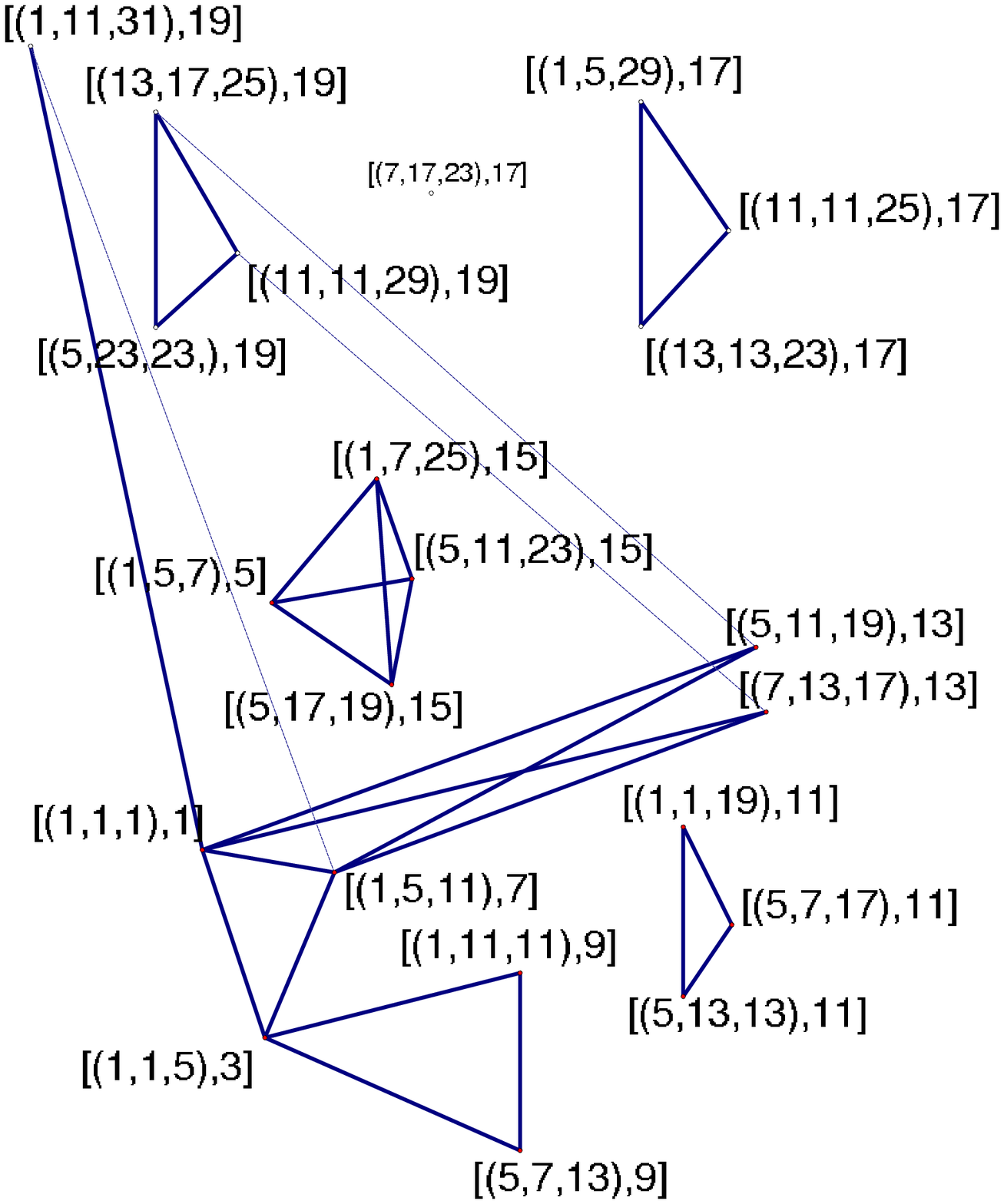,height=4.5in,width=3.5in}
\vspace{.1in}\\
{\small Figure 2: The graph ${\cal RT}$ , $d\le 19$.} \vspace{.1in}
\end{center}

A few questions related to this graph appear naturally at this
point. Is it connected? Is there a different characterization of the
existence of an edge between two vertices in terms of only $d_1$ and
$d_2$?  We do not have an answer to the second question. This graph
seems to have a fractal structure.

Each edge in this graph, determined by $[(a_1,b_1,c_1),d_1]$ and
$[(a_2,b_2,c_2),d_2]$, gives rise to a minimal tetrahedra (the side
lengths  are at most $max\{d_1,d_2\}\sqrt{2}$) which is determined
up to the set of isometric transformations that are generated by the
symmetries of the cube $\cal C(m)$ where m is the size of the
smallest ``cube" $\{0,1,\cdots,m\}^3$ containing the tetrahedron or
a translation of it.

\section{ Some preliminaries}

We would like to have a good estimate of the primitive solutions of
(\ref{plane}) which satisfy in addition $0<a\le b\le c$. Let us
observe that we cannot have $a=b=c$ unless $d=1$. So, the counting
in (\ref{2007hirschhorn}) via (\ref{legendresymbol}) would give what
we want if we can count the number of positive primitive solutions
of the following equation in terms of $d$:

\begin{equation}\label{x=y}
2a^2+c^2=3d^2.
\end{equation}

\n A similar description to the Pythagorean triples, which gives the
nature of the solutions of (\ref{x=y}), is stated next.

\begin{theorem}\label{characterizationofx=y}
For every positive integers $l$ and $k$ such that, gcd$(k,l)=1$ and
$k$ is odd, then $a$, $c$ and $d$ given by

\begin{equation}\label{parx=y}
d=2l^2+k^2 \ and \ \begin{cases} a=|2l^2+ 2kl-k^2|, \ c=|k^2+
4kl-2l^2|,\ if \ k\not\equiv l \ {\rm
(mod\ 3)}\\ \\
a=|2l^2-2kl-k^2|, \ c=|k^2-4kl-2l^2|,\ if \ k\not\equiv -l\  {\rm
(mod\ 3)}
\end{cases}
\end{equation}

\n constitute a positive primitive solution for {\rm (\ref{x=y})}.

Conversely, with the exception of the trivial solution $a=c=d=1$,
every positive primitive solution for {\rm (\ref{x=y})} appears in
the way described above for some $l$ and $k$.
\end{theorem}

\proof.\   First, one can check that (\ref{parx=y}) satisfy {\rm
(\ref{x=y})} for every $l$ and $k$. As a result it follows that $a$,
$c$ and $d$ are positive integers. Let $p$ be a prime dividing $a$,
$c$ and $d$. Then $p$ must divide $\pm a-d=2k(\pm l-k)$ and so $p$
is equal to $2$, $p$ divides $k$ or it divides $\pm l-k$. If $p=2$
then, $p$ must divide $k$ but this contradicts the assumption that
$k$ is odd.

In case $p$ is not equal to $2$ and it divides $k$, we see $p$ must
divide $l^2=(d-k^2)/2$. Since we assumed $gcd(l,k)=1$ it remains
that $p$ must divide $\pm l-k$. By our assumptions on $k$ and $l$,
$p$ cannot be equal to $3$. Then $p$ divides $\pm a+(\pm
l-k)^2=3l^2$. Because $p\not =3$ then  $p$ must divide $l^2$ and so
$p$ should divide $l$ and then $k$. This contradiction shows that
$a$, $c$ and $d$ cannot have prime common factors. So, we have a
primitive solution in (\ref{parx=y}).

For the converse, let us assume that $a$, $c$ and $d$ is a positive
primitive solution of (\ref{x=y}), which is different of the trivial
one. We denote by $u=\frac{a}{d}$ and $v=\frac{c}{d}$. Then the
point of rational coordinates $(u,v)$ (different of $(1,1)$)) is on
the ellipse $\frac{x^2}{3/2}+\frac{y^2}{3}=1$ (Figure 3) in the
first quadrant. This ellipse contains the following four points with
integer coordinates: $(1,1)$, $(-1,1)$, $(-1,-1)$ and $(1,1)$. This
gives the lines $y+1=t_1(x+1)$, $y+1=t_2(x-1)$, $y-1=t_3(x+1)$, and
$y-1=t_4(x-1)$, passing through $(u,v)$ and one of the points
mentioned above. Hence, the slopes $t_1$, $t_2$, $t_3$, and $t_4$,
are rational numbers. This gives expressions for the point $(u,v)$
in terms of $t_i$ ($i=1,...,4$). Let us assume that
$t_i=\frac{k_i}{l_i}$ with $k_i,l_i\in \mathbb Z$, written in the
reduced form. Then we must have

$$u=\frac{|2\pm 2t_i-t_i^2|}{2+t_i^2}=\frac{|2l_i^2\pm 2k_il_i-k_i^2|}{2l_i^2+k_i^2},
\ \ v=\frac{|t_i^2\pm 4t_i-2|}{2+t_i^2}=\frac{|k_i^2\pm
4k_il_i-2l_i^2|}{2l_i^2+k_i^2},$$

\begin{center}\label{figure2}
\epsfig{file=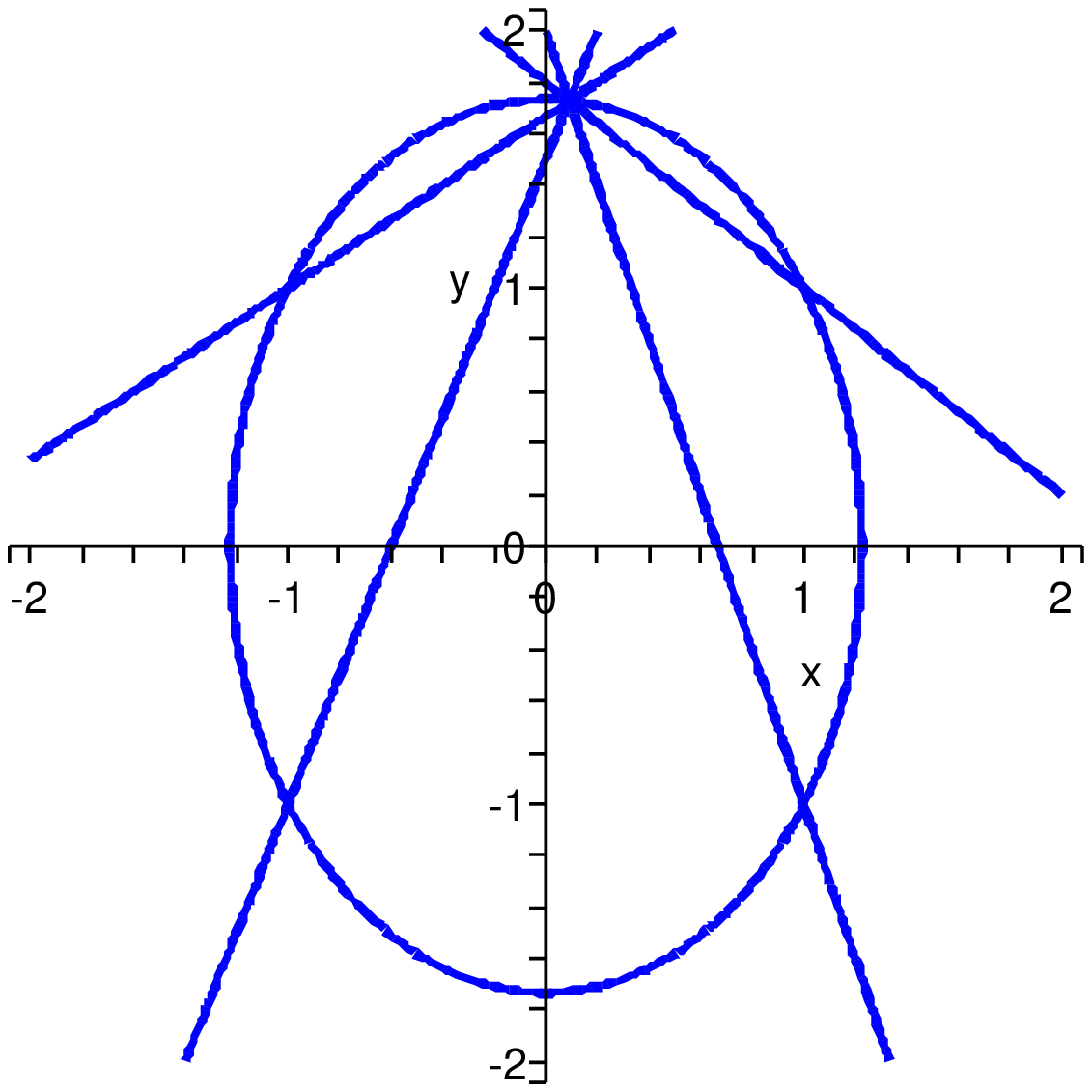,height=3in,width=3in}
\vspace{.1in}\\
{\small Figure 3: The ellipse $\frac{x^2}{3/2}+\frac{y^2}{3}=1$}
\vspace{.1in}
\end{center}

\n and so, these equalities give

\begin{equation}\label{representations}
\frac{a}{d}=\frac{|2l_i^2\pm 2k_il_i-k_i^2|}{2l_i^2+k_i^2},\ and\
\frac{c}{d}=\frac{|k_i^2\pm 4k_il_i-2l_i^2|}{2l_i^2+k_i^2},\
i=1,...,4.
\end{equation}

We claim that the function $t_i\to 2l_i^2+k_i^2$ ($i=1,...,4$) is
injective. If for some $2l_i^2+k_i^2=2l_j^2+k_j^2$ ($i\not=j$), that
would imply that the corresponding numerators in
(\ref{representations}) are equal. This gives enough infirmation to
conclude a contradiction. There are ${4 \choose 2}=6$ possibilities
here but we are going to include the details only in the case $i=1$
and $j=2$. The rest of the cases can be done in a similar fashion.
For this situation we have,
$2l_1^2+2k_1l_1-k_1^2=k_2^2+2k_2l_2-2l_2^2$ and
$k_1^2+4k_1l_1-2l_1^2=k_2^2-4k_2l_2-2l_2^2$. The first equality
implies
$$2k_1l_1=k_1^2+k_2^2+2k_2l_2-2l_1^2-2l_2^2=2k_1^2+2k_2l_2-4l_2^2$$ which
substituted into the second equality gives

$$6k_1^2=2k_2^2-8k_2l_2+8l_2^2\Leftrightarrow 3k_1^2=(k_2-2l_2)^2.$$

Because $\sqrt{3}$ is irrational, the last equality is impossible
for $k_1$, $k_2$, $l_2$ integers and $k_1$ nonzero. For the other
cases one will get a contradiction based on the facts that
$\sqrt{\frac{3}{2}}$ and $\sqrt{2}$ are irrational numbers.

A similar argument to the one in the first part of the proof shows
that the fractions in the right-hand side of the equalities of
(\ref{representations}) can be simplified only by a factor of $2$,
$3$ or $6$. Having four distinct possibilities in
(\ref{representations}) for the denominators, exactly one of the
fractions (simultaneously in the first and second equalities) must
be in reduced form. This one will give the wanted representation.
\eproof

Similar to Fermat's theorem about the representation of primes as a
sum of two squares and the number of such representations one can
show the next result.

\begin{theorem}\label{prime}{\bf (Fermat \cite{cox})}
An odd  prime $p$ can be written as $2x^2+y^2$ with $x,y\in \mathbb
Z$ if and only if $p\equiv 1$ or $3$ {\rm (mod\ 8)}. If $d=2^k\prod
p_i^{\alpha_i}\prod q_j^{\beta_j}$ is the prime factorization of
$d$, with $q_j$ primes as before and $p_i$ the rest of them, then
the number of representations $d=2x^2+y^2$ with $x,y\in \mathbb Z$
is either zero if not all $\alpha_i$ are even and otherwise given by
\begin{equation}\label{numberofrepr}
 \lfloor \frac{1}{2}\prod (\beta_i+1)\rfloor.
\end{equation}
The number of positive primitive representations $d=2x^2+y^2$ for
$d$ odd, i.e. $x,y\in \mathbb N$ and {\rm gcd$(x,y)=1$},  is equal
to

\begin{equation}\label{numberofrep2x2plusy2}
\Gamma_2(d)=\begin{cases} 0 \ \text{ if d is divisible by a prime
factor of the form 8s+5 or
8s+7}, \ s\ge 0,\\ \\
2^{k-1}\ \ \begin{cases}\text{ where k is the number of
distinct prime factors of d}   \\ \\
\text {of d of the form 8s+1, or  8s+3}\  (s\ge 0). \end{cases}
\end{cases}
\end{equation}

\end{theorem}

Putting the two results together (Theorem~\ref{prime} and
Theorem~\ref{primitivesolch}) we obtain the following proposition:

\begin{proposition}\label{counting}
For every $d$ odd, the number of representations of (\ref{plane})
which satisfy $0<a\le b\le c$ and gcd$(a,b,c)=1$ is equal to
\begin{equation}\label{numberofrepr}
\pi\epsilon(d)=\frac{\Lambda(d)+24\Gamma_2(3d^2)}{48}.
\end{equation}
\end{proposition}

A regular tetrahedron whose vertices are integers is said to be {\it
irreducible} if it cannot be obtained by an integer dilation and a
translation from a smaller one also with integer coordinates. An
important question at this point about irreducible tetrahedra is
included next.

{\it \begin{quote}
   Does every irreducible tetrahedron with integer coordinates have a
face with
    a normal vector $(a,b,c)$  satisfying $a^2+b^2+c^2=3d^2$, such that $d$ gives the
    side lengths $\ell$ of the tetrahedron by the formula
    $\ell=d\sqrt{2}$?  In other words, is there a face for which $k=1$ in the
Theorem~\ref{main}?
\end{quote}}

 Unfortunately the answer to this question is no.
  The following points together
with the origin, [-6677, -2672, 1445], [-5940, 4143, -1167], [-3837,
2595, 5688] form a regular tetrahedron of side-lengths equal to
$5187\sqrt{2}$ and the highest $d$ for the faces is 1729.
 We observe that $3$, $7$, $13$ and $19$ are the first three
distinct primes numbers of the form $u^2+3v^2$, $u,v\in \mathbb Z$.

\section{The Code}

The program is written in Maple code and it is based on the
Theorem~\ref{main}. The main idea is to create a list of
\underline{irreducible} regular tetrahedra that can be used to
generate all the others in $\{0,1,2,...,n\}^3$ by certain
transformations generating a partition for the set of all the
tetrahedra. Each such irreducible tetrahedron is constructed out of
the equation of one face using Theorem~\ref{main}. One important
problem that appears here is to make sure this list contains
distinct elements, elements which may appear theoretically
 in this list from four different constructions, one for each face.
It turns out that there is a simple way of making sure that this
doesn't happen.

Let us observe that if $gcd(m,n)=d>1$ then all the coordinates of
the vertices of the initial face are multiple of $d$ and by the
formula (\ref{fourthpoint}), so are the coordinates of the fourth
point of the tetrahedron. This is because the numbers $k$ in
(\ref{fourthpoint}) satisfy a Diophantine equation of the form
$k^2=m^2-mn+n^2$. We go one step further, if $k^2$ is of the form

\begin{equation}\label{mn}
m^2-mn+n^2=\left(\frac{m+n}{2}\right)^2+3\left(\frac{m-n}{2}\right)^2,
\end{equation}

\n then one can see that for $k$ even, the $gcd(m,n)\ge 2$. Also, if
$k$ is odd then both $m$ and $n$ must be odd and so we have a
representation of $k^2$ as $u^2+3v^2$, $u,v\in \mathbb Z$. If $k$ is
divisible by $3$ then it is easy to see that $3$ divides $u$ an $v$,
and this attracts $gcd(m,n)\ge 3$. Hence we are going to look only
for those odd values $k\le n$, which are not multiples of $3$, in
the Theorem~\ref{main}. This means that only one choice of signs in
(\ref{fourthpoint}) is useful. A similar fact to Theorem~\ref{prime}
takes place.

\begin{theorem}\label{prime2}{\bf (Fermat \cite{cox})}
A prime $p$ can be written as $x^2+3y^2$ with $x,y\in \mathbb Z$ if
and only if $p=3$ or $p\equiv 1$ {\rm (mod\ 3)}. If $d=\prod
p_i^{\alpha_i}\prod q_j^{\beta_j}$ with $q_j$ primes as before and
$p_i$ the rest of them, then the number of representations
$d=x^2+3y^2$ with $x,y\in \mathbb Z$ is either zero if not all
$\alpha_i$ are even and otherwise given by

\begin{equation}\label{numberofrepr3}
 \lfloor \frac{1}{2}\prod (\beta_i+1)\rfloor.
\end{equation}
The number of positive primitive representations $d=x^2+3y^2$ for
$d$ odd, i.e. $x,y\in \mathbb N$ and {\rm gcd$(x,y)=1$},  is equal
to

\begin{equation}\label{numberofrepx2plus3y2}
\Gamma_3(d)=\begin{cases} 0 \ \text{ if d is even or divisible by a
prime factor of the form 3s+2}, \ s\ge 0,\\ \\
2^{k-1}\ \ \begin{cases}\text{ where k is the number of
distinct prime factors of d such as 3 or}   \\ \\
\text {of the form 3s+1}\  (s\ge 2). \end{cases}
\end{cases}
\end{equation}

\end{theorem}

As a result of these facts we first calculate all $k\le n$ such that
$k^2=m^2-mn+n^2$ has a solution with $gcd(m,n)=1$.

\begin{quote}
\n {\it kvalues}:=proc(n)

\n \ \ local i,j,k,L,a,p,q,r,m,mm;

\n \ \ L:={};mm:=floor((n+1)/2);

\n \ \ for i from 2 to mm do a:=ifactors(2i-1); k:=nops(a[2]);r:=0;

\n \ \ \ for j from 1 to k do

\n \ \ \  \    m:=a[2][j][1]; p:=m mod 3;

\n \ \ \  \    if m=3 then r:=1 fi;

\n \ \ \  \    if r=0 and p=2 then r:=1 fi;

 \n \ \ \   od;

\n \ \ \ if r=0 then L:=L union $\{2i-1\}$;fi;

\n \ \ od;

\n L:=L union $\{1\}$; L:=convert(L,list);

\n end:
\end{quote}

\n The result of this procedure for $n=100$ is [1, 7, 13, 19, 31,
37, 43, 49, 61, 67, 73, 79, 91, 97]. These are all the natural
numbers less than 100 which are primes of the form $3s+1$ or
products of such primes. We noticed that it takes only a fraction of
a second to execute this procedure if we limit $n$ to be less than
10000, although it may time consuming for big numbers. So an
alternative solution to this procedure may use a similar result to
that in Theorem~\ref{characterizationofx=y}, to describe all the
solutions of (\ref{mn}). For each $k$ found by the previous
procedure, there are usually at least eighteen solutions of
(\ref{mn}) if signs and order are counted, but if we impose the
conditions $gcd(m,n)=1$, $0<m, n$  and $2m<n$, we slice these
solutions by a factor of $18$. Such a solution is going to be
referred to as a {\it primitive} solution and these primitive
solution of (\ref{mn}) can be calculated with the following
procedure.

\begin{quote}

{\it listofmn}:=proc(k)

\n  \ local a,b,i,x,m,n,nx,L;

\n \ x:=[isolve($k^2=m^2-mn+n^2$)]; nx:=nops(x); L:=\{ \};

\n \ \ \ for i from 1 to nx do

\n \ \ \ \ \ if lhs(x[i][1])=m then a:=rhs(x[i][1]);
b:=rhs(x[i][2]);

\n \ \ \ \ \ \ \ \ else b:=rhs(x[i][1]); a:=rhs(x[i][2]); fi;

\n \ \ \ \ \ if gcd(a,b)=1 and $a>0$ and $b>0$ and $2a<b$ then L:=L
union {[a,b]};fi;

\n \ \ od;L;

\n end:

\end{quote}
For example, if $k=91$ we get two primitive solutions, $m=1991$,
$n=9095$,  and $m=3401$, $n=9440$. It turns out that it is enough to
know just the primitive solutions of (\ref{mn}) in order to find all
integer solutions. For each $k\le n$, output of the procedure {\em
kvalues}, we need to find the values of $d$ as in Theorem\ref{main},
which are only restricted to two conditions: $d$ must be an odd
positive integer and $dk\le n$. The last restriction follows from
the fact that the sides of the regular tetrahedron constructed from
say $d$ and $k$, with a primitive solution $(m,n)$ of (\ref{mn}),
must be equal to $dk\sqrt{2}$. One can see that the biggest regular
tetrahedron inscribed in the cube $[0,n]^3$ has sides equal to $
n\sqrt{2}$ (see Proposition~2.1 in \cite{ejic}). The next procedure
is then very simple.

\begin{center}\label{figure4}
\epsfig{file=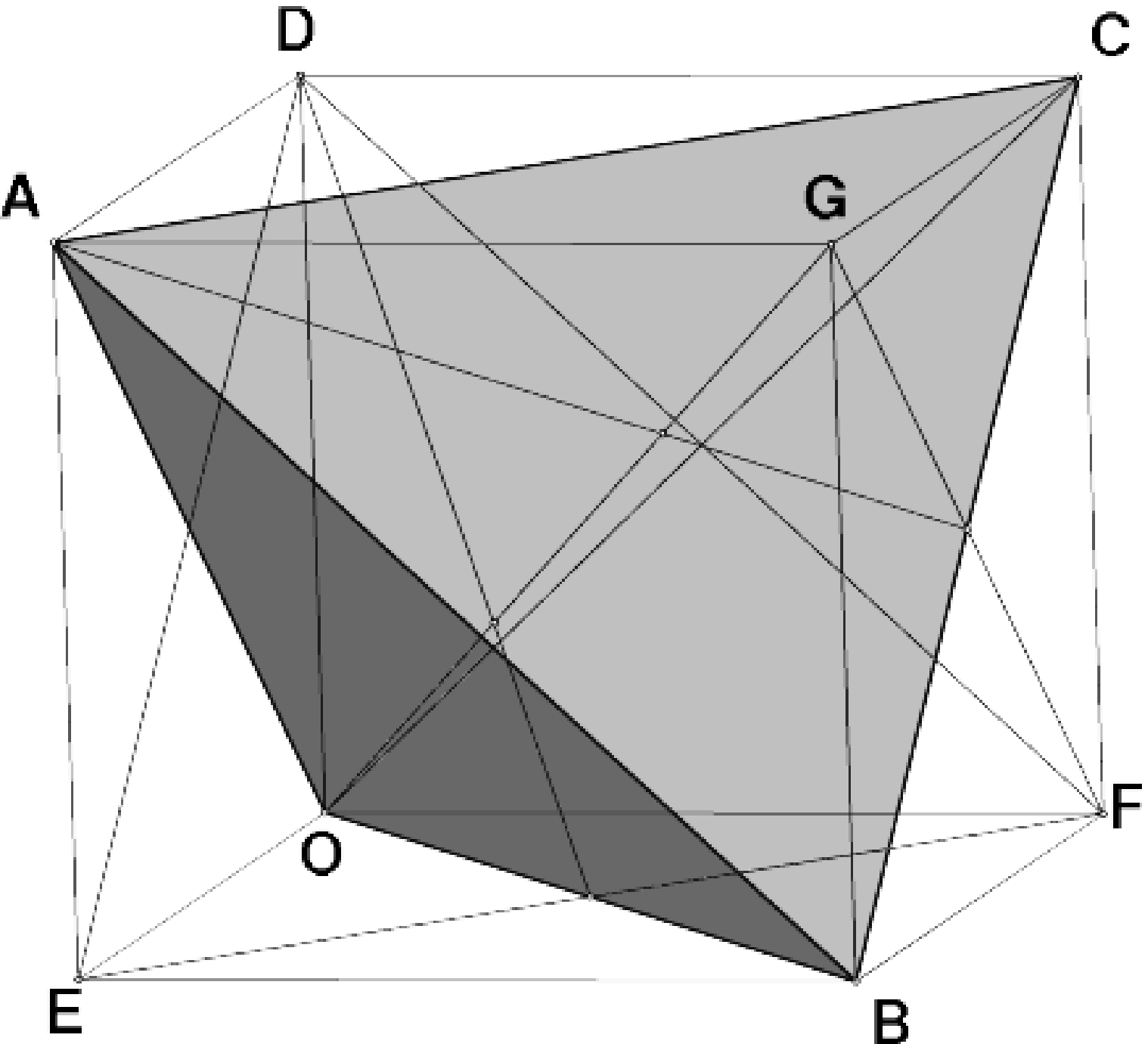,height=2in,width=2in}
\vspace{.1in}\\
{\small Figure 4: Largest tetrahedrons in a cube} \vspace{.1in}
\end{center}

\begin{quote}
\n {\it determined}:=proc(n)\

\n  local i,x,m,L,j;

\n x:=kvalues(n);m:=nops(x);

\n for i from 1 to m do

 \n \ \ j:=1;L[i]:=$\{  \}$;

\n  \ \  while (2j-1)x[i]$<$=n  do

\n  \ \  L[i]:=L[i] union $\{2j-1\}$;j:=j+1; od;

\n od; [seq(L[i],i=1..m)];

\n end:

\end{quote}

\n Next, we need to find all primitive solutions $(a,b,c)$ of the
equation $a^2+b^2+c^2=3d^2$. The number of such solution is given by
(\ref{numberofrepr}).

\begin{quote}
{\it abcsol}:=proc(d)

\n local i,j,k,m,u,x,y,sol,cd;

\n sol:=$\{ \}$;

\n for i from 1 to d do

\n \ \  u:=[isolve($3d^2-i^2=x^2+y^2$)];k:=nops(u);

\n  \ \ \  for j from 1 to k do

\n  \ \ \ \     if rhs(u[j][1])$>$=i and rhs(u[j][2])$>$=i then

\n   \ \ \ \ \ \  \ \     cd:=gcd(gcd(i,rhs(u[j][1])),rhs(u[j][2]));

\n  \ \ \ \ \ \ \  if cd=1 then sol:=sol union
$\{sort([i,rhs(u[j][1]),rhs(u[j][2])])\}$;fi;

\n  \ \   \ \   fi;

\n od; od;

\n convert(sol,list);

\n end:

\end{quote}

\n For example, if $d=2009$ we get 294 solutions, as seen before, of
which one of them is $a=1$, $b=1159$ and $c=3281$. Next, we are
going to use the construction of an equilateral triangle in the
plane of equation $ax+by+cz=0$ with $a^2+b^2+c^2=3d^2$ using the
formulae (\ref{paramtwo}) with $m$ and $n$ given by the procedure
{\em listofmn}. The fourth point is then calculated using the
formula (\ref{fourthpoint}). We are using only two of the possible
values of $m$ and $n$  ($m'=m$, $n'=n$ and $m'=n$,  $n'=n-m$, in
order to obtain two equilateral triangles that share a side,
$\overrightarrow{OQ}$) since all other tetrahedra as in the figure
below, can be obtained from these two by a simple translation, and
as a result they will translate into the same minimal tetrahedron
inside the first quadrant (see \cite{ejirt}).

\begin{center}\label{figure5}
\epsfig{file=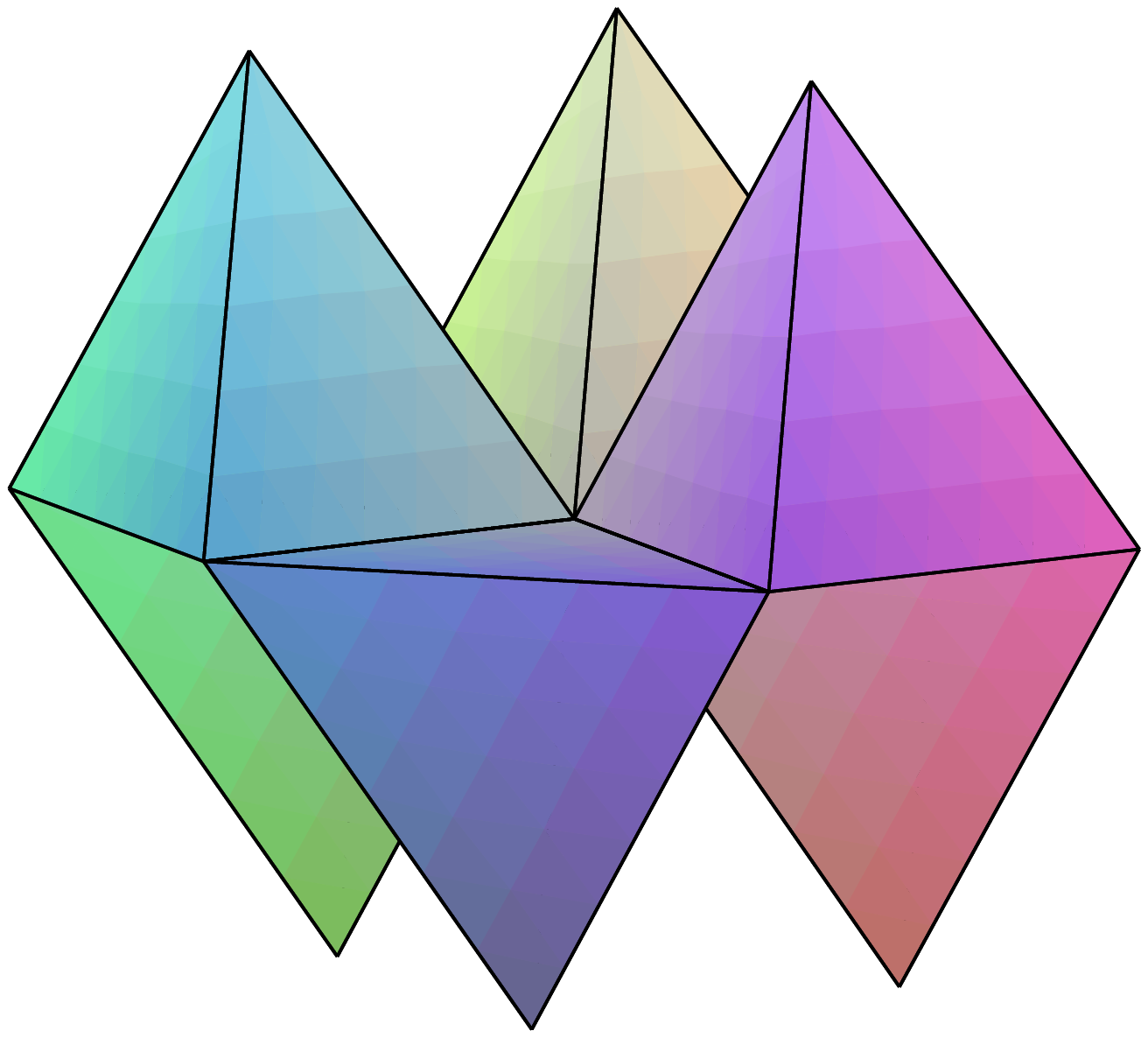,height=2in,width=2in}
\vspace{.1in}\\
{\small Figure 5: All six tetrahedrons with one face in ${\cal
P}_{a,b,c}$} \vspace{.1in}
\end{center}

\begin{quote}

\n {\it findpar}:=proc(a,b,c,mm,nn)

\n local
i,j,sol,mx,nx,r,s,my,ny,q,d,u,v,w,x,y,z,mu,nu,mv,nv,ef,ns,mz,nz,

\n mw,nw,om1,om2,l,uu,t,R1,R2,fc,k;\par

\n $q:=a^2+b^2$; $k:=sqrt(mm^2-mm*nn+nn^2)$;

\n sol:=convert({isolve($2*q=x^2+3*y^2)$},list); ns:=nops(sol);
d:=sqrt($(a^2+b^2+c^2)/3$); ef:=0;

\n for i from 1 to ns do

\n \ if ef=0 then r:=rhs(sol[i][1]); s:=rhs(sol[i][2]);

\n \ \    uu:=$(s^2+3*r^2-2*q)^2$; {\em  if $uu>0$ then
t:=s;s:=r;r:=t; fi;}

\n \ \ mx:=-$(d*b*(3*r+s)+a*c*(r-s))/(2*q)$;nx:=-$(r*a*c+d*b*s)/q$;

\n \ \   my:=$(d*a*(3*r+s)-b*c*(r-s))/(2*q)$;ny:=-$(r*b*c-d*a*s)/q$;

\n \ \    mz:=(r-s)/2;nz:=r;
mu:=nx;mv:=ny;mw:=nz;nu:=nx-mx;nv:=ny-my; nw:=nz-mz;

\n \ \  if mx=floor(mx) and nx=floor(nx) and my=floor(my) and
ny=floor(ny) then

\n\ \  u:=mu*m-nu*n;v:=mv*m-nv*n;w:=mw*m-nw*n; \par\n
\par\n \ \ x:=mx*m-nx*n;y:=my*m-ny*n;z:=mz*m-nz*n;
\par\n \ \ R1:=[(x+u-2*a*k)/3,(v+y-2*b*k)/3,(z+w-2*c*k)/3];
\par\n \ \ R2:=[(x+u+2*a*k)/3,(v+y+2*b*k)/3,(z+w+2*c*k)/3];
\par\n \ \ fc:=subs(m=mm,n=nn,R1[1]); if fc=floor(fc) then
\par\n \ \ om1:=subs(m=mm,n=nn,[[u,v,w],[x,y,z],R1]); else
\par\n \ \ om1:=subs(m=mm,n=nn,[[u,v,w],[x,y,z],R2]);fi;
\par\n \ \ fc:=subs(m=nn,n=nn-mm,R1[1]); if fc=floor(fc) then
\par\n \ \ om2:=subs(m=nn,n=nn-mm,[[u,v,w],[x,y,z],R1]); else
\par\n \ \ om2:=subs(m=nn,n=nn-mm,[[u,v,w],[x,y,z],R2]);fi; ef:=1; fi;fi; od;
\par\n om1,om2; end:

\end{quote}

Before we translate these two tetrahedrons we need a small
sub-routine for subtraction of  two vectors.

\begin{quote} {\it subtrv}:=proc(U,V) local W;
\par\n \ \ W[1]:=U[1]-V[1];W[2]:=U[2]-V[2];W[3]:=U[3]-V[3]; [W[1],W[2],W[3]];
end:
\end{quote}

The next procedure translates a tetrahedron which comes as output of
{\it findpar} into the positive octant of the space in such a way
that for each component at least one of the tetrahedron's vertex has
a zero coordinate on that component. Let us observe that this
operation is invariant to translations of the tetrahedron; this
justifies the choice of looking at only two tetrahedrons out of six
in the procedure {\it findpar}.

\begin{quote} {\it
tmttopq}:=proc(T) local i,a,b,c,v,O,TR;
\par\n \ \  a:=min(T[1][1],T[2][1],T[3][1],0);
\par\n \ \  b:=min(T[1][2],T[2][2],T[3][2],0);
\par\n \ \  c:=min(T[1][3],T[2][3],T[3][3],0); O:=[0,0,0];
v:=[a,b,c];
\par\n \ \  TR:=\{subtrv(O,v),subtrv(T[1],v),subtrv(T[2],v),subtrv(T[3],v)\};
end:
\end{quote}

Next, we calculate the size of smallest cube $C_m:=[0,m]^3$ which
contains the tetrahedron resulted from the {\it tmttopq}.

\begin{quote} {\it
mscofmt}:=proc(Q) local a,b,c,T,m; T:=convert(Q,list);
\par\n \ \  a:=max(T[1][1],T[2][1],T[3][1],T[4][1]);
\par\n \ \  b:=max(T[1][2],T[2][2],T[3][2],T[4][2]);
\par\n \ \  c:=max(T[1][3],T[2][3],T[3][3],T[4][3]); m:=max(a,b,c); end:
\end{quote}

The tetrahedron obtained as a result of  {\it tmttopq} is then
transformed within the cube found above through all the
translations, rotations and symmetries of the cube. We denote this
orbit of $T$, by $O(T)$.

\begin{quote} {\it
orbit1}:=proc(T) local
i,k,T1,a,b,c,x,T2,T3,T4,T5,T6,T7,T8,T9,T10,T11,T12,T13,T14,T15,T16,
\par\n \ \ T17,T18,T19,T20,T21,T22,T23,T24,S,Q,d,a1,b1,c1; Q:=convert(T,list);
\par\n \ \  d:=mscofmt(T); T1:=T; T2:=\{seq([Q[k][2],Q[k][3],Q[k][1]],k=1..4)\};
\par\n \ \  T3:=\{seq([Q[k][1],Q[k][3],Q[k][2]],k=1..4)\};
\par\n \ \  T4:=\{seq([Q[k][1],Q[k][2],d-Q[k][3]],k=1..4)\};
\par\n \ \ T5:=\{seq([Q[k][2],Q[k][3],d-Q[k][1]],k=1..4)\};
\par\n \ \ T6:=\{seq([Q[k][1],Q[k][3],d-Q[k][2]],k=1..4)\};
\par\n \ \ T7:=\{seq([Q[k][1],d-Q[k][2],Q[k][3]],k=1..4)\};
\par\n \ \ T8:=\{seq([Q[k][2],d-Q[k][3],Q[k][1]],k=1..4)\};
\par\n \ \ T9:=\{seq([Q[k][1],d-Q[k][3],Q[k][2]],k=1..4)\};
\par\n \ \ T10:=\{seq([d-Q[k][1],Q[k][2],Q[k][3]],k=1..4)\};
\par\n \ \ T11:=\{seq([d-Q[k][2],Q[k][3],Q[k][1]],k=1..4)\};
\par\n \ \ T12:=\{seq([d-Q[k][1],Q[k][3],Q[k][2]],k=1..4)\};
\par\n \ \ T13:=\{seq([Q[k][1],d-Q[k][2],d-Q[k][3]],k=1..4)\};
\par\n \ \ T14:=\{seq([Q[k][2],d-Q[k][3],d-Q[k][1]],k=1..4)\};
\par\n \ \ T15:=\{seq([Q[k][1],d-Q[k][3],d-Q[k][2]],k=1..4)\};
\par\n \ \ T16:=\{seq([d-Q[k][1],d-Q[k][2],Q[k][3]],k=1..4)\};
\par\n \ \ T17:=\{seq([d-Q[k][2],d-Q[k][3],Q[k][1]],k=1..4)\};
\par\n \ \ T18:=\{seq([d-Q[k][1],d-Q[k][3],Q[k][2]],k=1..4)\};
\par\n \ \ T19:=\{seq([d-Q[k][1],Q[k][2],d-Q[k][3]],k=1..4)\};
\par\n \ \ T20:=\{seq([d-Q[k][2],Q[k][3],d-Q[k][1]],k=1..4)\};
\par\n \ \ T21:=\{seq([d-Q[k][1],Q[k][3],d-Q[k][2]],k=1..4)\};
\par\n \ \ T22:=\{seq([d-Q[k][1],d-Q[k][2],d-Q[k][3]],k=1..4)\};
\par\n \ \ T23:=\{seq([d-Q[k][2],d-Q[k][3],d-Q[k][1]],k=1..4)\};
\par\n \ \ T24:=\{seq([d-Q[k][1],d-Q[k][3],d-Q[k][2]],k=1..4)\};
\par\n \ \ S:=\{T1,T2,T3,T4,T5,T6,T7,T8,T9,T10,T11,T12,T13,T14,T15,T16,T17,
\par\n \ \ T18,T19,T20,T21,T22,T23,T24\};S; end:
\end{quote}

\begin{quote} {\it
orbit}:=proc(T) local S,Q,T1; Q:=convert(T,list);
\par\n \ \ T1:=\{[Q[1][3],Q[1][2],Q[1][1]],[Q[2][3],Q[2][2],Q[2][1]],
\par\n \ \ [Q[3][3],Q[3][2],Q[3][1]],[Q[4][3],Q[4][2],Q[4][1]]\};
\par\n \ \ S:=orbit1(T) union orbit1(T1); S; end:
\end{quote}

We recall from \cite{ejic} a few variables that we are going to use
in this calculation also. The theorem used there applies as well to
this case because it is  a pure set theoretic result. The meaning of
those variables here is:
\begin{enumerate}
  \item $n$ -the dimension of the cube $C_n=[0,n]^3$,
  \item $m$ -the maximum of all the coordinates in a tetrahedron $T$ computed by {\it
tmttopq},
  \item $\alpha(T)$ -the cardinality of $O(T)$ within $C_m$,
  \item $\beta(T)$ -the cardinality of $O(T)\cap [O(T)+e_1]$,
  \item $\gamma(T)$ -the cardinality of $[O(T)+e_1] \cap [O(T)+e_2] $.
\end{enumerate}

\begin{theorem}(Theorem 2.2 in\cite{ejic}) \label{calcofall}
The number $f(T,n)$ of all tetrahedrons that can be obtained from
$T$ within a cube $C_n$ by translations, rotations, or symmetries,
is given by
\begin{equation}\label{eq2}
f(T,n)=(n+1-m)^3\alpha-3(n+1-m)^2(n-m)\beta+3(n+1-m)(n-m)^2\gamma ,
\end{equation}
 for
all $n\ge m$.
\end{theorem}

Hence, we need to calculate $\alpha$, $\beta$ and $\gamma$.

\begin{quote} {\it  addvec}:=proc(U,V)
 local W;
\par\n \ \  W[1]:=U[1]+V[1];W[2]:=U[2]+V[2];W[3]:=U[3]+V[3];[W[1],W[2],W[3]];
\par\n \ \  end:
\end{quote}

\begin{quote} {\it  addvect}:=proc(T,v)
\par\n \ \  local i,Q;Q:=\{  \};
\par\n \ \  for i from 1 to 4 do
\par\n \ \ Q:=Q union {addvec(T[i],v)};
\par\n \ \ od;Q;end:
\end{quote}

\begin{quote} {\it transl}:=proc(T)
\par\n \ \ local S,Q,i,j,k,a2,b2,c2,a,b,c,d;
\par\n \ \ Q:=convert(T,list);
\par\n \ \ a:=max(Q[1][1],Q[2][1],Q[3][1],Q[4][1]);
\par\n \ \  b:=max(Q[1][2],Q[2][2],Q[3][2],Q[4][2]);
\par\n \ \ c:=max(Q[1][3],Q[2][3],Q[3][3],Q[4][3]);
\par\n \ \ d:=max(a,b,c);
\par\n \ \ a2:=d-a;b2:=d-b;c2:=d-c;
\par\n \ \ S:=orbit(T);
\par\n \ \ for i from 0 to a2 do
\par\n \ \ \  for j from 0 to b2 do
\par\n \ \  \ \   for k from 0 to c2 do
\par\n \ \      \ \ \  S:=S union orbit(addvect(T,[i,j,k]));
\par\n \ \      od; od; od; S; end:
\end{quote}
This last procedure gives the value of $\alpha$. Then $\beta$ is
calculated by the following.

\begin{quote} {\it
intersalongE1}:=proc(T)
\par\n \ \  local S,m,i,S1,S2;
\par\n \ \  S2:=transl(T);S:=convert(S2,list);m:=nops(S);S1:=\{ \};
\par\n \ \ for i from 1 to
m do S1:=S1 union {addvect(S[i],[0,0,1])};
\par\n \ \ od; S2 intersect S1; end:
\end{quote}

Then the variable $\gamma$ is given by the procedure:

\begin{quote} {\it
intersalongE2}:=proc(T) local S,m,i,S1,S2,S3,S4;
\par\n \ \   S2:=transl(T);S:=convert(S2,list);m:=nops(S);S1:=\{ \};
\par\n \ \ \  for i from 1 to
m do S1:=S1 union \{addvect(S[i],[0,0,1])\};
\par\n \ \ \ od; S3:=\{\};
\par\n \ \ \ \  for i
from 1 to m do S3:=S3 union \{addvect(S[i],[0,1,0])\}; od; S1
intersect S3; end:
\end{quote}

The function in Theorem~\ref{calcofall} is then implemented by

\begin{quote}
 {\it f}:=(n,d,alpha,beta,gamma)$\rightarrow (n-d+1)^3*alpha-3*(n-d)*(-d+1+n)^2*beta+3*gamma*(n-d+1)*(n-d)^2$:
\end{quote}

\begin{quote}
 {\it notetraincn}:=proc(T,n)
\par\n \ \ local d,x,y,z,w;
\par\n \ \ d:=mscofmt(T);
\par\n \ \ x:=nops(transl(T));y:=nops(intersalongE1(T));w:=nops(intersalongE2(T));
\par\n \ \ if $n\ge d$ then
\par\n \ \ z:=f(n,d,x,y,w);else z:=0;
\par\n \ \ fi;z;end:
\end{quote}

Finally we put together a list of irreducible tetrahedrons whose
orbits under the operations above are pairwise disjoint. The next
four procedures are pretty simple and one can figure out what they
do. They are used in the code of {\it
 ExtendList}.

\begin{quote}
{\it  distance}:=proc(A,B) local C;
\par\n \ \  C:=subtrv(A,B); sqrt($C[1]^2+C[2]^2+C[3]^2$);
\par\n end:
\end{quote}
\begin{quote}
{\it checkrt}:=proc(T) local d1,d2,d3,d4,d5,d6,Q,D;
\par\n \ d1:=distance(T[4],T[1])/sqrt(2); d2:=distance(T[4],T[2])/sqrt(2);
\par\n \ d3:=distance(T[4],T[3])/sqrt(2); d4:=distance(T[1],T[2])/sqrt(2);
\par\n \ d5:=distance(T[1],T[3])/sqrt(2); d6:=distance(T[2],T[3])/sqrt(2);
\par\n min(d1,d2,d3,d4,d5,d6); end:
\end{quote}

\begin{quote}
 {\it crossproductt}:=proc(U,V) local x,i,j,k,d,i1,j1,k1,d1;
\par\n \ \ i:=U[2]*V[3]-U[3]*V[2];j:=U[3]*V[1]-U[1]*V[3];k:=U[1]*V[2]-V[1]*U[2];
\par\n \ \ d:=gcd(i,j);d1:=gcd(d,k);i1:=i/d1;j1:=j/d1;k1:=k/d1;
sqrt(($i1^2+j1^2+k1^2)/3$); end:
\end{quote}

\begin{quote}
{\it facesnew}:=proc(T) local N1,N2,N3,N4,U,V,W;
\par\n \ \ U:=subtrv(T[1],T[2]);V:=subtrv(T[1],T[3]);W:=subtrv(T[1],T[4]);
\par\n \ \ N1:=crossproductt(U,V);N2:=crossproductt(V,W);N3:=crossproductt(U,W);
\par\n \ \ N4:=crossproductt(subtrv(U,V),subtrv(U,W)); max(N1,N2,N3,N4); end:
\end{quote}

\begin{quote}
{\it
 ExtendList}:=proc(n,L,mm,nn)
\par\n local i,sol,nsol,alpha,beta,gammma,nel,mt,ttp1,ttp2,Or1,Or,tnel,NL,intcard,cio,normals,length,k;
\par\n  nel:=nops(L);
\par\n k:=sqrt($mm^2-mm*nn+nn^2$);
\par\n  sol:=abcsol(n);nsol:=nops(sol);
\par\n  tnel:=nel;Or:=\{  \};
\par\n  \  for i from 1 to nsol do
\par\n  \ \   mt:=findpar(sol[i][1],sol[i][2],sol[i][3],mm,nn); ttp1:=tmttopq(mt[1]); ttp2:=tmttopq(mt[2]);
\par\n  \ \  normals:=facesnew(ttp1);length:=checkrt(ttp1);
\par\n \ \      if k=1 or $normals<length$ then
\par\n    cio:=evalb(ttp1 in Or);
\par\n      if cio=false then
\par\n       Or:=Or union transl(ttp1); Or1:=transl(ttp2);intcard:=nops(Or intersect Or1);
\par\n      if $intcard>0$ then NL[tnel+1]:=[n,mscofmt(ttp1),ttp1,sol[i]];
\par\n       tnel:=tnel+1;else
\par\n          NL[tnel+1]:=[n,mscofmt(ttp1),ttp1,sol[i]];
\par\n         NL[tnel+2]:=[n,mscofmt(ttp2),ttp2,sol[i]];
\par\n       tnel:=tnel+2;
\par\n       fi;
\par\n       fi;
\par\n   \ fi;
\par\n  \ normals:=facesnew(ttp2);length:=checkrt(ttp2);
\par\n  \ if k=1 or $normals<length$ then
\par\n   \ \  ttp2:=tmttopq(mt[2]);
\par\n   \ \   cio:=evalb(ttp2 in Or);
\par\n   \ \ \    if cio=false then
\par\n  \ \ \      Or:=Or union transl(ttp2); Or1:=transl(ttp1);intcard:=nops(Or intersect Or1);
\par\n  \ \ \      if $intcard>0$ then NL[tnel+1]:=[n,mscofmt(ttp1),ttp1,sol[i]];
\par\n  \ \ \ \      tnel:=tnel+1;else
\par\n  \ \ \ \        NL[tnel+1]:=[n,mscofmt(ttp1),ttp1,sol[i]];
\par\n  \ \ \ \        NL[tnel+2]:=[n,mscofmt(ttp2),ttp2,sol[i]];
\par\n   \ \      tnel:=tnel+2;
\par\n   \ \   fi;
\par\n  \     fi;
\par\n   fi;
\par\n  od;
\par\n  [seq(L[i],i=1..nel),seq(NL[j],j=nel+1..tnel)];
\par\n  end:
\end{quote}

Once this list is computed for every $n$ we can add up all the
contributions for each tetrahedron in the list. A few procedures are
used in the code of the main calculation called {\it calculation}.

\begin{quote}
{\it multbyfactorv}:=proc(v,k)
\par\n \  local w;
 \par\n \ \ \ w[1]:=v[1]*k;w[2]:=v[2]*k;w[3]:=v[3]*k;
\par\n \ \ \ [w[1],w[2],w[3]];
\par\n  end:
\end{quote}

\begin{quote}
{\it  multbyfactor}:=proc(T,k)
\par\n  \ local i,NT,Q;NT:=\{\};
\par\n \ \ Q:=convert(T,list);
\par\n  \ \ \ for i from 1 to 4 do
\par\n \ \ \ \ NT:=NT union {multbyfactorv(Q[i],k)};
\par\n \ \ \ od;NT;
\par\n end:
\end{quote}

\begin{quote}
{\it  addup}:=proc(n,L)
\par\n local i,j,m,mm,nt,k,T,Q,alpha,beta,gammma,d;nt:=0;
\par\n \ m:=nops(L);k:=floor((n-1)/2); i:=1;
\par\n  \ \ while $i\le m$  do
\par\n \ \  if $L[i][2]\le n$ then
\par\n  \ \ \     mm:=floor(n/L[i][1]);T:=L[i][3];j:=1;
\par\n  \ \ \ \    while $j\le mm$  and $L[i][2]*j\le n$ do
\par\n \ \ \  \   Q:=multbyfactor(T,j);d:=L[i][2]*j;
\par\n \ \ \ \     alpha:=nops(transl(Q));beta:=nops(intersalongE1(Q));
\par\n \ \ \  \ gammma:=nops(intersalongE2(Q));nt:=nt+f(n,d,alpha,beta,gammma);
\par\n \ \    j:=j+1;
\par\n  \  od;fi;
\par\n \   i:=i+1;
\par\n od;
\par\n nt/2;
\par\n end:
\end{quote}

In the above procedure, {\it addup}, we get every irreducible from
the list $L$ together with all their appropriate dilations to
compute their contribution in the cube $[0,n]^3$ using the formula
given by Theorem~\ref{calcofall}. The result is divided by two to
match the sequence introduced in $A103158$. Finally, the list $L$ is
calculated in terms of $n$ in two steps. First, construct the
irreducible tetrahedrons using values of of $mm=0$ and $nn=1$ in
{\it ExtendList} and then we take care of the other possible values
of $mm$ and $nn$, making sure that we do not have duplicates by
looking at the biggest $d$ that shows up from each face.

\begin{quote}
\n {\it calculation}:=proc(n) local i,j,j1,j2,k,k1,L,x,xx,y,B,ii,T;
\par\n  x:=kvalues(n);y:=determined(n); xx:=floor((n+1)/2);
\par\n  \ i:=nops(x);L:=[];L:=ExtendList(1,L,0,1);
\par\n \ \ for j from 1 to xx do
\par\n \ \ \  L:=ExtendList(2*j+1,L,0,1);
\par\n  \ \ \  \ od; for j from 1 to i do
\par\n  \ \ \ \   k:=listofmn(x[j]);k1:=nops(k);
\par\n  \  \ \ \ \ for j2 from 1 to k1  do
 \par\n \ \ \ \ \   for j1 from 1 to nops(y[i]) do
 \par\n  \ \ \   L:=ExtendList(y[j][j1],L,k[j2][1],k[j2][2]);
 \par\n  \ \ od;
\par\n \ od;
\par\n od; for ii from 1 to n do B[ii]:=addup(ii,L); print([ii,B[ii]]); od;
\par\n T:=[seq(B[ii],ii=1..n)]; T; end:
\end{quote}

The result of the {\it calculation(100)} gives in less than a few
hours of computations:
                             [1, 1]
                             [2, 9]
                            [3, 36]
                            [4, 104]
                            [5, 257]
                            [6, 549]
                           [7, 1058]
                           [8, 1896]
                           [9, 3199]
                           [10, 5145]
                           [11, 7926]
                          [12, 11768]
                          [13, 16967]
                          [14, 23859]
                          [15, 32846]
                          [16, 44378]
                          [17, 58977]
                          [18, 77215]
                          [19, 99684]
                          [20, 126994]
                          [21, 159963]
                          [22, 199443]
                          [23, 246304]
                          [24, 301702]
                          [25, 366729]
                          [26, 442587]
                          [27, 530508]
                          [28, 631820]
                          [29, 748121]
                          [30, 880941]
                         [31, 1031930]
                         [32, 1202984]
                         [33, 1395927]
                         [34, 1612655]
                         [35, 1855676]
                         [36, 2127122]
                         [37, 2429577]
                         [38, 2765531]
                         [39, 3137480]
                         [40, 3548434]
                         [41, 4001071]
                         [42, 4498685]
                         [43, 5044606]
                         [44, 5641892]
                         [45, 6294195]
                         [46, 7005191]
                         [47, 7778912]
                         [48, 8620242]
                         [49, 9533105]
                         [50, 10521999]
                         [51, 11591474]
                         [52, 12746562]
                         [53, 13992107]
                         [54, 15332971]
                         [55, 16775590]
                         [56, 18324372]
                         [57, 19985523]
                         [58, 21765013]
                         [59, 23668266]
                         [60, 25702480]
                         [61, 27873699]
                         [62, 30188259]
                         [63, 32655348]
                         [64, 35281418]
                         [65, 38074085]
                         [66, 41040495]
                         [67, 44188592]
                         [68, 47525856]
                         [69, 51061295]
                         [70, 54804647]
                         [71, 58763604]
                         [72, 62949850]
                         [73, 67371219]
                         [74, 72037311]
                         [75, 76958126]
                         [76, 82143618]
                         [77, 87606245]
                         [78, 93355379]
                         [79, 99403446]
                        [80, 105762770]
                        [81, 112443331]
                        [82, 119456581]
                        [83, 126814970]
                        [84, 134532746]
                        [85, 142621185]
                        [86, 151093691]
                        [87, 159964136]
                        [88, 169245226]
                        [89, 178954039]
                        [90, 189102295]
                        [91, 199706864]
                        [92, 210781424]
                        [93, 222341631]
                        [94, 234402515]
                        [95, 246978962]
                        [96, 260093046]
                        [97, 273757925]
                        [98, 287989943]
                        [99, 302809940]
                        [100, 318235290]

We observe a similar behavior with the sequence
$\frac{ln(ET(n))}{\ln(n+1)}$, in \cite{ejic}:

\begin{center}\label{figure6}
\epsfig{file=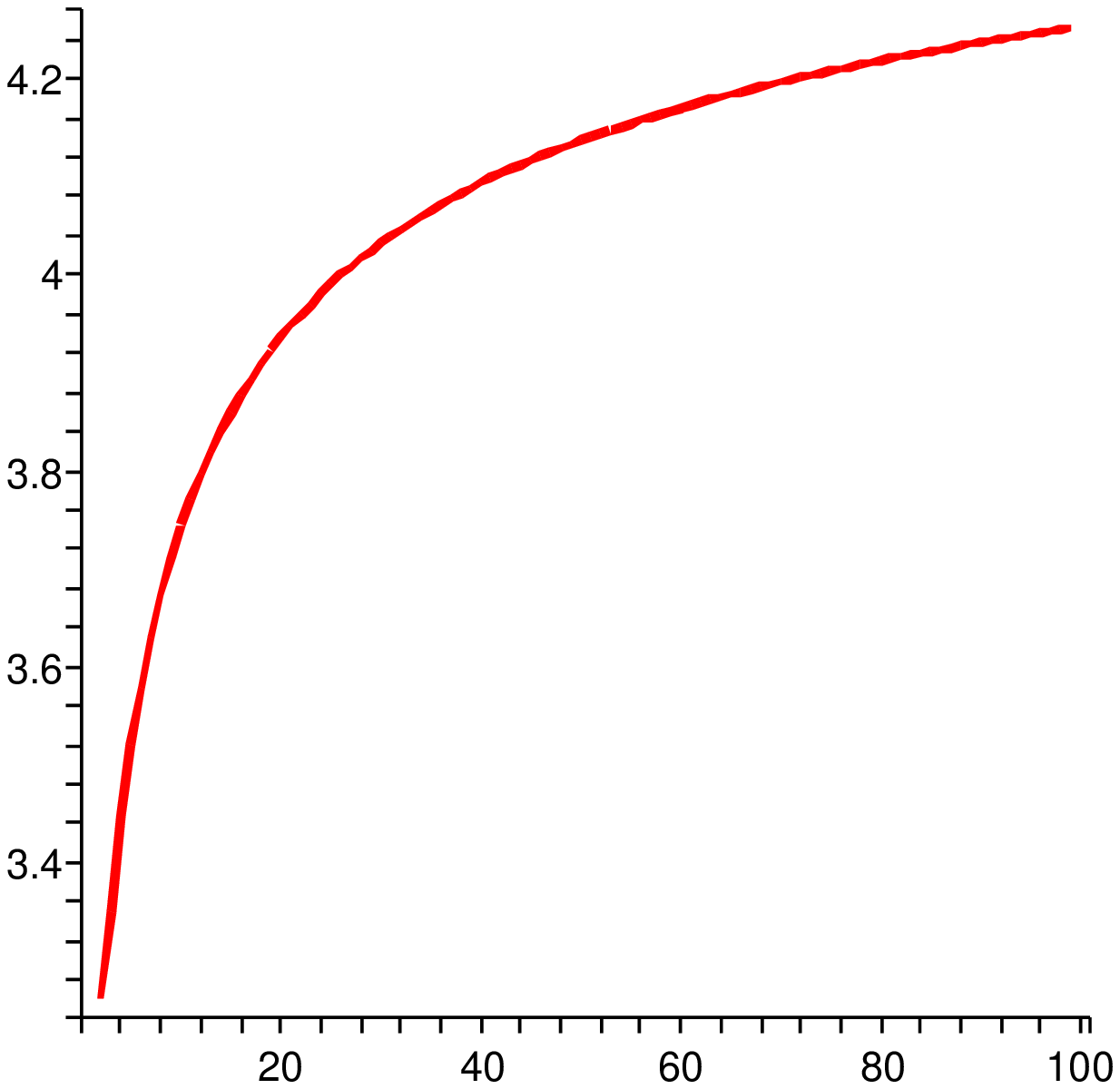,height=3.5in,width=3.5in}
\vspace{.1in}\\
{\small Figure 6: The graph $\frac{ln(T(n)/2)}{\ln(n+1)}$ , $1\le
n\le 100$.} \vspace{.1in}
\end{center}

\end{document}